\documentclass[11pt]{article}
\usepackage[left=1in, right=1in, top=1in, bottom=1in]{geometry}
\usepackage{lineno}
\usepackage{enumitem}
\usepackage{authblk}
\usepackage{graphicx}	
\usepackage{amsmath, amsthm, amssymb}
\usepackage{xspace}
\usepackage{comment}
\usepackage{hyperref}
\usepackage{xcolor}
\usepackage{float}

\usepackage[labelsep=period]{caption}
\makeatletter
\def\widebreve{\mathpalette\wide@breve}
\def\wide@breve#1#2{\sbox\z@{$#1#2$}%
     \mathop{\vbox{\m@th\ialign{##\crcr
\kern0.08em\brevefill#1{0.8\wd\z@}\crcr\noalign{\nointerlineskip}%
                    $\hss#1#2\hss$\crcr}}}\limits}
\def\brevefill#1#2{$\m@th\sbox\tw@{$#1($}%
  \hss\resizebox{#2}{\wd\tw@}{\rotatebox[origin=c]{90}{\upshape(}}\hss$}
\makeatletter

\theoremstyle{plain}
\newtheorem{theorem}{Theorem}[section]
\newtheorem{proposition}[theorem]{Proposition}
\newtheorem{corollary}[theorem]{Corollary}
\newtheorem{lemma}[theorem]{Lemma}
\newtheorem{observation}[theorem]{Observation}

\theoremstyle{definition}
\newtheorem{definition}[theorem]{Definition}
\newtheorem{remark}[theorem]{Remark}
\newtheorem{problem}[theorem]{Problem}
\newtheorem{question}[theorem]{Question}

\newcommand{\Lem}[1]{Lemma~\ref{#1}\xspace}

\newcommand{\Prop}[1]{Proposition~\ref{#1}\xspace}
\newcommand{\Thm}[1]{Theorem~\ref{#1}\xspace}
\newcommand{\Def}[1]{Definition~\ref{#1}\xspace}

\DeclareMathOperator{\Spec}{Spec}

\title{Directed Ramsey and Anti-Ramsey Schemes and the Flexible Atom Conjecture
\footnote{ML was partially supported by J. A. Grochow's NSF award CISE-2047756. We wish to thank J.N. Cooper for helpful discussions, as well as the anonymous referee for their very careful reading and thoughtful suggestions which have greatly improved the exposition of our work.} }

\author[1]{Jeremy F. Alm}
\author[2,3]{Michael Levet\footnote{Corresponding author: levetm@cofc.edu.}}
\affil[1]{Department of Mathematics, Lamar University}
\affil[2]{Department of Computer Science, College of Charleston}
\affil[3]{Department of Computer Science, University of Colorado Boulder}

\begin{document}
\maketitle

\begin{abstract}
In this paper, we shed new light on the Flexible Atom Conjecture. We first give finite representation results for relation algebras $33_{37}, 35_{37}$, $77_{83}$, $78_{83}$, $80_{83}$, $82_{83}$, $83_{83}$, $1310_{1316}$, $1313_{1316}$, $1315_{1316}$, and $1316_{1316}$. Prior to our paper, only $83_{83}$ and  $1316_{1316}$ were known to be finitely representable. We accomplish this by generalizing the notion of a relation algebra generated by a Ramsey scheme to the directed (antisymmetric) setting, and then showing that each of these algebras embeds into a finite directed (anti-)Ramsey scheme. The notion of a directed (anti-)Ramsey scheme may be of independent interest.

We complement our upper bounds with some lower bounds. Namely, we show that any square representation of $31_{37}$ requires at least $16$ points, any square representation of $33_{37}$ requires at least $11$ points, and any square representation of $35_{37}$ requires at least $12$ points. Our technique adapts previous work of Alm, et. al. (\textit{Algebra Universalis} 2022), in that we examine the combinatorial structure induced by the flexible atom.
\end{abstract}

\textbf{Kewords.} Flexible Atom Conjecture, Directed Ramsey schemes, Representations. \\

\textbf{Mathematics Subject Classification.} 03G15.

\thispagestyle{empty}

\newpage

\setcounter{page}{1}

\section{Introduction}
The \textit{Flexible Atom Conjecture} states that every finite integral relation algebra with a flexible atom (i.e., an atom that does not participate in any forbidden diversity cycles) is representable over a finite set. Jipsen, Maddux, \& Tuza showed that the finite symmetric integral relation algebras in which every diversity atom is flexible, are finitely representable. In particular, the algebra with $n$ flexible symmetric atoms is representable over a set of size $(2+o(1))n^2$ \cite{Jipsen}. We note that if all cycles are present, then all diversity atoms are flexible. Hence, the case considered in \cite{Jipsen} is intuitively the \textit{big end} of the Flexible Atom Conjecture.

The other extreme is when just enough cycles are present for one atom to be flexible. This case was handled by Alm, Maddux, \& Manske, who exhibited a representation of $32_{65}$ over a set of size $416,714,805,914$ \cite{AMM}. Dodd \& Hirsch \cite{DH} subsequently improved the bound the representation size to $63,432,274,896$. This was subsequently improved to 8192 by J.F. Alm \& D. Sexton (unpublished), and later $3432$ by \cite{AlmAndrews}. Finally, in \cite{alm2021improved}, the authors exhibited a representation over a set of size $1024$, as well as the first polynomial upper bound on a more general family of relation algebras $A_{n}$ obtained by splitting the non-flexible diversity atom of $6_{7}$ into $n$ symmetric atoms.

There are strikingly few lower bounds in the literature. Jipsen, Maddux, \& Tuza exhibited a lower bound of $n^{2} + n + 1$ for the relation algebra $\mathfrak{E}_{n+1}(1,2,3)$ \cite{Jipsen}. In \cite{alm2021improved}, the authors showed that any representation of $A_{n}$ requires at least $2n^{2} + 4n + 1$ points, which is asymptotically double the trivial lower bound of $n^{2} + 2n + 3$. The key technique involved analyzing the combinatorial substructure induced  by the flexible atom. In the special case of $A_{2} = 32_{65}$, the authors used a SAT solver to further improve the lower bound. Namely, they showed that $32_{65}$ is not representable on a set of fewer than $26$ points.

The notion of a Ramsey scheme was first introduced (but not so named) by Comer \cite{ComerMonochrome}, where he used them to obtain finite representations of relation algebras. Kowalski \cite{KowalskiRamsey} later introduced the term \textit{Ramsey relation algebra} to refer to the (abstract) relation algebras obtained from embeddings into Ramsey schemes. With the sole exception of an alternate construction of the $3$-color algebra using $(\mathbb{Z}/4\mathbb{Z})^{2}$ \cite{WHITEHEAD1975399}, all known constructions have been due to the guess-and-check finite field method of Comer \cite{ComerMonochrome}. Intuitively, Comer's method takes a positive integer $m$ and considers finite fields $\mathbb{F}_{q}$, where $q \equiv 1 \pmod{2m}$. We next consider the multiplicative subgroup $H \leq \mathbb{F}_{q}^{\times}$ of order $(q-1)/m$, and check whether the cosets of $\mathbb{F}_{q}^{\times}/H$ yield a representation (namely, taking the cosets to be the  atoms of the relation algebra). The $m$-color Ramsey number 
\[
R_{m}(3) = R(\underbrace{3, \ldots, 3}_{m})
\]

\noindent provides an upper bound on $q$, restricting the search space for the finite fields to be considered. Comer was able to construct $m$-color Ramsey schemes for all $1 \leq m \leq 5$  \cite{ComerMonochrome}. In 2011, Maddux leveraged improvements in computing power to produce constructions for $m = 6, 7$ \cite{MadduxAMS2011}. The cases of $m = 8, 13$ were elusive. Alm \& Manske \cite{AlmManskeRamsey} produced constructions over prime fields for all $m \leq 400$, with the exception of $m = 8, 13$. However, they ruled out $m = 8$ by considering all primes up through the Ramsey bound $R_{8}(3)$. Independently, Kowalski produced constructions for all $m \leq 120$, except for $m = 8, 13$. Kowalski's constructions included non-prime fields. Additionally, Kowalski also ruled out $m = 8$ by checking all primes up through the Ramsey bound $R_{8}(3)$ \cite{KowalskiRamsey}. In subsequent work, Alm gave constructions for $401 \leq m \leq 2000$. He also substantially improved the upper bound on $p$ with respect to $m$, finally showing that no construction over prime fields exists for $m = 13$ \cite{Alm2017401AB}.

Ramsey schemes also have close connections with other combinatorial structures such as association schemes, coherent configurations, and permutation groups \cite{ComerMonochrome, ComerCombinatorial}. 

\noindent \\ \textbf{Main Results.} In this paper, we rule out several possible counter-examples to the Flexible Atom Conjecture. Precisely, we show that several families of relation algebras are finitely representable.

\begin{theorem} \label{ThmMainRepresentation}
The following relation algebras are finitely representable:
\begin{enumerate}[label=(\alph*)]
\item $33_{37}$ over a set of size $29$.
\item $35_{37}$ over a set of size $3221$.
\item $77_{83}$ over a set of size $29$.
\item $78_{83}$ over a set of size $33791$.
\item $80_{83}$ over a set of size $67.$
\item $82_{83}$ over a set of size $3221$.
\item $83_{83}$ over a set of size $37$. 
\item $1310_{1316}$ over a set of size $67$.
\item $1313_{1316}$ over a set of size $3221$.
\item $1315_{1316}$ over a set of size $33791$. 
\item $1316_{1316}$  over a set of size $73$.
\end{enumerate}
\end{theorem}

Prior to our paper, none of the following relation algebras were known to be finitely representable: $33_{37}, 35_{37}$, $77_{83}$, $78_{83}$, $80_{83}$, $82_{83}$, $1310_{1316}$, $1313_{1316}$, and $1315_{1316}$. In particular, this rules out each of these relation algebras as a counter-example to the Flexible Atom Conjecture. As every diversity atom of $83_{83}$  and $1316_{1316}$ are flexible, both of these algebras are finitely representable by a slight generalization of \cite{Jipsen}.

\begin{remark}
\Thm{ThmMainRepresentation} can be written more compactly by observing the following:
\begin{itemize}
    \item $33_{37} \leq 77_{83}$, where $77_{83}$ is the $2$-color directed anti-Ramsey algebra.

    \item $35_{37} \leq 82_{83}$, which are subalgebras of the $10$-color directed Ramsey algebra.

    \item $1313_{16}$ is also a subalgebra of the $10$-color directed Ramsey algebra.

    \item $80_{83}$ and $1310_{16}$ are both subalgebras of the $3$-color directed anti-Ramsey algebra.

    \item $78_{83}$ admits a representation by splitting the symmetric atoms of the $31$-color symmetric Ramsey algebra, into asymmetric pairs.

    \item $83_{83}$ embeds into the $4$-color directed anti-Ramsey algebra.

    \item $1315_{1316}$ is a subalgebra of another type of directed Ramsey algebra; see Section~\ref{1315} for details.

    \item $1316_{1316}$ is a subalgebra of another type of directed Ramsey algebra; see Section~\ref{1316} for details.
\end{itemize}

We emphasize the representation sizes in the statement of \Thm{ThmMainRepresentation}, as our approach yields small representations in all of these cases. In particular, we conjecture that the representations given for $33_{37}$ and $77_{83}$ are indeed optimal. In general, finite representations need not be small- for instance, the first finite representation of $32_{65}$ in \cite{AMM} was over a set of size $416,714,805,914$; it took $14$ years and a number of follow-up papers \cite{DH, AlmAndrews, alm2021improved} to finally obtain a representation on $1024$ points.
\end{remark}

We prove Theorem \ref{ThmMainRepresentation} in two steps. First, we generalize the notions of Ramsey relation algebras to the directed (antisymmetric) case. This construction may be of independent interest. Next, we show that most of the relation algebras in Theorem \ref{ThmMainRepresentation} embed into some representable directed (anti-)Ramsey algebra. This second step provides our finite representations.  

Next, we complement our upper bounds with some lower bounds for $31_{37}, 33_{37}$, and $35_{37}$. Each of these relation algebras forbids certain blue triangles. We refer to the intransitive triangle $(x, y), (y, z), (z, x)$ as a \textit{graph-theoretic cycle}, and the transitive triangle $(x, y), (y, z), (x, z)$ as a \textit{graph-theoretic chain}. As we will be formulating our results in the language of graph theory, we will refer to these objects as the \textit{cycle} and \textit{chain} throughout. We note that Kramer \& Maddux \cite{KramerMaddux} refer to the graph-theoretic cycle as a \textit{loop}.\footnote{We will later define a family $\text{Cy}_{n}$ of relation algebras that forbid the graph-theoretic cycle. To avoid confusion with Lyndon algebras (so that we may use $\mathrm{Cy}_{n}$ rather than $\mathrm{L}_{n}$), we will make reference to the graph-theoretic cycle rather than adopting the term \textit{loop} from Kramer \& Maddux \cite{KramerMaddux}.}

Note that the graph-theoretic cycle is an instance of a relation-algebraic diversity cycle (see Definition~\ref{def:DiversityCycle}) $bb\breve{b}$. The relation algebra $31_{37}$ forbids both blue cycles and blue chains, while $33_{37}$ forbids only blue cycles and $35_{37}$ forbids only blue chains. Precisely, we show the following.

\begin{theorem} \label{thm:MainLowerBounds}
\noindent
\begin{enumerate}[label=(\alph*)]
\item If $n < 16$, then $n \not \in \text{Spec}(31_{37})$.
\item If $n < 11$, then $n \not \in \text{Spec}(33_{37})$.
\item If $n < 12$, then $n \not \in \text{Spec}(35_{37}).$
\end{enumerate}
\end{theorem}

We prove Theorem \ref{thm:MainLowerBounds} by adapting the techniques of \cite{alm2021improved}. Namely, we  examine the structure induced by the flexible atom, which provides initial lower bounds.

We also obtain lower bounds for two different generalizations of $31_{37}, 33_{37},$ and $35_{37}$. The first way generalizes these relation algebras in the direction of Ramsey schemes, in which we forbid the respective triangles that all have the same shade of blue. Let $\rho$ be our representation. We say that the chain $(x, y), (y, z), (x, z)$ is monochromatic blue if $(x, y), (y, z), (x, z) \in \rho(b_{i})$. The monochromatic cycle is defined analogously. Define the relation algebras $\mathrm{MonoCy}_{n}, \mathrm{MonoCh}_{n}, \mathrm{MonoCych}_{n}$ to have atoms $1^{\prime}, a, b_{1}, \breve{b_{1}}, \ldots, b_{n}, \breve{b_{n}}$, where $\mathrm{MonoCy}_{n}$ forbids the monochromatic blue cycles, $\mathrm{MonoCh}_{n}$ forbids the monochromatic blue chains, and $\mathrm{MonoCy}_{n}$ forbids the monochromatic blue cycles, and $\mathrm{MonoCych}_{n}$ forbids both monochromatic blue cycles and monochromatic blue chains. Observe that $\mathrm{MonoCych}_{1} = 31_{37}$, $\mathrm{MonoCy}_{1} = 33_{37}$, and $\mathrm{MonoCh}_{1} = 35_{37}$.

We may alternatively generalize $31_{37}, 33_{37},$ and $35_{37}$ in the direction of \cite{AMM} by forbidding the respective blue triangles, where the shades of blue need not all be the same. We consider three families of relation algberas, which we label as $\mathrm{Cy}_n$, $\mathrm{Ch}_n$, $\mathrm{Cych}_n$, with the atoms $1^{\prime}, a, b_{1}, \breve{b_{1}}, \ldots, b_{n}, \breve{b_{n}}$. Now $\mathrm{Cy}_{n}$ forbids blue cycles, $\mathrm{Ch}_{n}$ forbids blue chains, and $\mathrm{Cych}_{n}$ forbids both blue cycles and blue chains. Thus, $\mathrm{Cych}_{1} = 31_{37}, \mathrm{Cy}_{1} = 33_{37}$, and $\mathrm{Ch}_{1} = 35_{37}$.

\begin{proposition} \label{PropLowerBounds}
Let $n \geq 2$, and let $\mathcal{C} \in \{ \mathrm{MonoCy}, \mathrm{MonoCh}, \mathrm{MonoCych}, \mathrm{Cy}, \mathrm{Ch}, \mathrm{Cych}\}$. Let $f_{\mathcal{C}}(n)$ be the minimum-sized square representation of $\mathcal{C}_{n}$. We have that $f(\mathcal{C}_{n}) \geq 4n^{2} + 4n + 3$. 
\end{proposition}

The proof of Proposition \ref{PropLowerBounds} relies on examining the combinatorial structure induced by the flexible atom (a red edge). 

Now let $\{u, v\}$ be a red edge (precisely, if $\rho$ is our representation and $(u,v), (v, u) \in \rho(a)$), and let $x_{1}, x_{2}$ be elements meeting different blue-blue needs for $\{u,v\}$. In the cases of $\mathrm{Cy}_{n}, \mathrm{Ch}_{n}$, and $\mathrm{Cych}_{n}$, $\{x_{1}, x_{2}\}$ must necessarily be colored red. Forbidding the polychromatic blue chain provides enough structure for us to leverage the techniques of \cite[Theorem~3.5]{alm2021improved} in order to analyze the needs of $\{x_{1}, x_{2}\}$. This analysis provides stronger lower bounds, which we think are of independent interest. Precisely, we obtain the following.

\begin{theorem}
Let $n \geq 2$. We have that:
\begin{enumerate}[label=(\alph*)]
\item $f_{\mathrm{Cych}}(n) \geq 8n^{2} + 8n +  3$.

\item $f_{\mathrm{Ch}}(n) \geq 6n^{2} + 8n + 5.$
\end{enumerate}
\end{theorem}

\noindent \\ \textbf{Bonus Result.} In the process of our work on directed (anti-)Ramsey schemes, we came across a way to improve the upper bound for multiplicative van der Waerden numbers. Alm \cite{AlmArithmeticProgressions} previously showed that $\text{VW}^{\times}(n) \leq (1+o(1))n^{4}$. We improve this upper bound to the following:

\begin{proposition}
$\text{VW}^{\times}(n) \leq (1-o(1))n^{4} + O(n^{7/2})$.
\end{proposition}

\section{Preliminaries}

\subsection{Relation Algebras}
We recall some preliminary notions of relation algebras.

\begin{definition} \label{Def:RAs}
A \textit{relation algebra} is an algebra $\langle A, \land, \vee, \neg, 0, 1, \circ, \breve{}, 1^{\prime} \rangle$ that satisfies the following.
\begin{itemize}
\item $\langle A, \land, \vee, \neg, 0, 1 \rangle$ is a Boolean algebra, with $\neg$ the unary negation operator, $0$ the identity for $\vee$, and $1$ the identity for $\land$.

\item $\langle A, \circ, 1^{\prime} \rangle$ is a monoid, with $1^{\prime}$ our nullary identity. That is, relational composition is associative, and there is an identity relation $1^{\prime}$. 

\item $\breve{}$ \, is the unary converse operation, which is an anti-involution with respect to composition. Namely, $\breve{\breve{a}} = a$ for all $a \in A$, and $\widebreve{a \circ b} = \breve{b} \circ \breve{a}$ for all $a, b \in A$.

\item Converse and composition both distribute over disjunction. Precisely, for all $a, b, c \in A$, we have that:
\begin{align*}
&\widebreve{(a \vee b)} = \breve{a} \vee \breve{b}, \text{ and } \\
&(a \vee b) \circ c = (a \circ c) \vee (b \circ c).
\end{align*}

\item (Triangle Symmetry) For all $a, b, c \in A$, we have that:
\[
(a \circ b)\land c = 0 \iff (\breve{c} \circ a)\land \breve{b} = 0 \iff (b\circ \breve{c})\land \breve{a} = 0.
\]

Note that the Triangle Symmetry axiom defines an equivalence relation on triples of diversity atoms that corresponds to the symmetries of the triangle. 

\end{itemize}
\end{definition}

\noindent When the relation algebra is understood, we simply write $A$ rather than $\langle A, \land, \vee, \neg, 0, 1, \circ, \breve{}, 1^{\prime} \rangle$.

\begin{definition}
Let $A$ be a relation algebra. We say that $a \in A$ is an \textit{atom} if $a \neq 0$ and $b < a \implies b = 0$. Furthermore, we say that $a$ is a \textit{diversity atom} if $a$ also satisfies $a \land 1^{\prime} = 0$. 
\end{definition}

\begin{definition} \label{def:DiversityCycle}
For diversity atoms $a, b, c$, the triple $(a,b, c)$-- usually denoted $abc$-- is called a (relation algebraic) \textit{diversity cycle}. We say that $abc$ is \textit{forbidden} if $(a \circ b) \land c = 0$ and \textit{mandatory} if $a \circ b \geq c$. We say that a diversity atom is \textit{symmetric} if $a = \breve{a}$. \end{definition}

\begin{remark}
We note that for integral relation algebras, the composition operation $\circ$ is determined by the diversity cycles. 
\end{remark}

\begin{definition}
Let $f$ be a diversity atom. We say that $f$ is \textit{flexible} if for all diversity atoms $a, b$, we have that $abf$ is mandatory.
\end{definition}

\begin{definition}
We say that a relation algebra $A$ is \textit{representable} if there exists a set $U$ and an equivalence relation $E \subseteq U \times U$ such that $A$ embeds into
\[
\langle 2^{E}, \cup, \cap, ^{c}, \circ, ^{-1}, E, \emptyset, \text{Id}_{U} \rangle.
\]

\noindent Here, $^{c}$ is set complementation, and $^{-1}$ is the relational inverse. 
\end{definition}

In this paper, we will only be concerned with simple relation algebras. So there exists a set $U$ such that $E = U \times U$. We call such a representation \textit{square}.

\begin{definition}
Let $A$ be a finite relation algebra. Denote:
\[
\text{Spec}(A) := \{ \alpha \leq \omega : A \text{ has a square representation over a set of cardinality } \alpha\}.
\]
\end{definition}

\begin{definition}
Let $n \geq 1$. We define the following infinite families $\mathrm{MonoCh}_{n}, \mathrm{MonoCy}_{n},$ $\mathrm{MonoCych}_{n}$, $\mathrm{Ch}_{n},$ $\mathrm{Cy}_{n}$, and $\mathrm{Cych}_{n}$ of integral relation algebras with the diversity atoms $a, b_{1}, \ldots, b_{n}, \breve{b_{1}}, \ldots, \breve{b_{n}}$ with $a$ symmetric. We now specify the forbidden diversity cycles. 
\begin{itemize}
    \item $\mathrm{MonoCh}_{n}$ forbids the diversity cycles are $b_ib_ib_i$ for all $1\leq i \leq n$ (as well as cycles equivalent by Triangle Symmetry -- see \Def{Def:RAs}).
    
    \item $\mathrm{MonoCy}_{n}$ forbids the diversity cycles are $b_ib_i \breve{b_i}$ for all $1\leq i \leq n$ (as well as cycles equivalent by Triangle Symmetry).
    
    \item $\mathrm{MonoCych}_{n}$ forbids the diversity cycles are $b_ib_i \breve{b_i}, b_{i}b_{i}b_{i}$ for all $1\leq i \leq n$ (as well as cycles equivalent by Triangle Symmetry).
    
    \item $\mathrm{Ch}_{n}$ forbids the diversity cycles are $b_ib_jb_k$ for all $1\leq i,j,k \leq n$ (as well as cycles equivalent by Triangle Symmetry).
    
    \item $\mathrm{Cy}_{n}$ forbids the diversity cycles are $b_ib_j \breve{b_k}$ for all $1\leq i,j,k \leq n$ (as well as cycles equivalent by Triangle Symmetry).
    
    \item $\mathrm{Cych}_{n}$ forbids the diversity cycles are $b_ib_j \breve{b_k}, b_{i}b_{j}b_{k}$ for all $1\leq i,j,k \leq n$ (as well as cycles equivalent by Triangle Symmetry).
\end{itemize}
\end{definition}

\begin{remark}
While deleting the symmetric diversity atom $a$ in the relation algebras from Definition 2.8 would better align with the notion of a Ramsey scheme, we include $a$ so as to create classes of algebras that fall under the umbrella of the Flexible Atom Conjecture.
\end{remark}

\subsection{Fourier Analysis}
In this section, we recall some preliminary notions of Fourier Analysis over finite Abelian groups. For a more thorough treatment, we refer to \cite{IrelandRosen, BabaiFourier}. 

Let $G$ be a finite Abelian group of order $n$, written additively. A \textit{multiplicative group character} of $G$ is a group homomorphism $\chi : G \to \mathbb{C}^{\times}$. Observe that for all $a \in G$, we have that $\chi(a)^{n} = \chi(na) = \chi(0) = 1.$ So the values of $\chi$ are $n$th roots of unity. In particular, we observe that $\chi(a)^{-1} = \chi(-a) = \overline{\chi(a)}$. The \textit{principal character} of $G$, denoted $\chi_{1}$, maps $\chi_{1}(a) = 1$ for all $a \in G$. In particular, for any non-principal character $\chi$, $\sum_{g \in G} \chi(g) = 0$. This yields the \textit{First Orthogonality Relation} for characters, which states that for characters $\chi, \psi$:
\[
\sum_{g \in G} \chi(g) \overline{\psi(g)} = \begin{cases}
n & : \text{if } \chi = \psi, \\ 0 & : \text{ otherwise.}
\end{cases}
\]

\noindent Now for any two characters $\chi, \psi$, we have that $(\chi \psi)(a) = \chi(a)\psi(a)$. Thus, the characters themselves form a group under pointwise multiplication. We refer to this group as the \textit{dual group} of $G$, denoted $\hat{G}$.

If $G$ is cyclic, we may write the characters of $G$ as $\chi_{j}(a) = \omega^{ja}$, where $\omega$ is a primitive $n$th root of unity and $0 \leq j < n$. In particular, the following properties hold:
\begin{enumerate}[label=(\alph*)]
\item $\chi_{i} = \chi_{j}$ if and only if $i \equiv j \pmod{n}$,
\item $\chi_{j} = \chi_{1}^{j}$,
\item $\widehat{\mathbb{Z}/n\mathbb{Z}} = \{ \chi_{0}, \ldots, \chi_{n-1}\}$,
\item $\mathbb{Z}/n\mathbb{Z} \cong \widehat{\mathbb{Z}/n\mathbb{Z}}$.
\end{enumerate}

\noindent \\ Now suppose that $G = H_{1} \oplus H_{2}$ is Abelian. If $\varphi_{i} : H_{i} \to \mathbb{C}^{\times}$ is a character of $H_{i}$, then $\varphi(h_{1}, h_{2}) = \varphi_{1}(h_{1}) \oplus \varphi_{2}(h_{2})$ is a character of $G$. Conversely, all such characters of $G$ are obtained in this way. Now if $G$ is an arbitrary finite Abelian group, then $G$ decomposes as a direct sum of cyclic groups. Thus, for any finite Abelian group $G$, we have that $G \cong \hat{G}$.

The \textit{character table} $C$ for a finite Abelian group $G$ is an $n \times n$ matrix, where $C_{ij} = \chi_{i}(a_{j})$. The matrix $A = \frac{1}{\sqrt{n}} C$ is a unitary matrix. Hence, we have the \textit{Second Orthogonality Relation} for characters, which states that for any $a, b \in G$:
\[
\sum_{\chi \in \hat{G}} \chi(a) \overline{\chi(b)} = \begin{cases} n & : \text{if } a = b, \\ 0 & : \text{otherwise}. \end{cases}
\]

\noindent Let $\mathbb{C}[G]$ denote the \textit{complex group algebra} of $G$, where
\[
\mathbb{C}[G] = \left \{ \sum_{g \in G} c_{g} g : c_{g} \in \mathbb{C} \right \}.
\]

\noindent It is well-known that $\mathbb{C}[G] \cong \mathbb{C}^{G}$, where $\mathbb{C}^{G}$ is the set (and in fact, vector space) of functions $f : G \to \mathbb{C}$. We note that $\hat{G}$ forms an orthonormal basis for $\mathbb{C}[G]$. The \textit{Fourier transform} is a function $\mathcal{F} : \mathbb{C}[G] \to \mathbb{C}[G]$, where for a function $f \in \mathbb{C}[G] \cong \mathbb{C}^{G}$, we map $\mathcal{F}(f) : f \mapsto \hat{f}$, where $\hat{f} : \hat{G} \to \mathbb{C}$ is given by:
\[
\hat{f}(\chi) = \sum_{a \in G} f(a)\chi(a).
\]

\noindent We may recover $f$ from $\hat{f}$ as follows:
\[
f(a) = \frac{1}{n} \cdot \sum_{\chi \in \hat{G}} \hat{f}(\chi) \overline{\chi(a)}.
\]

\section{Directed Ramsey Schemes}  \label{sec:RamseySchemes}

In this section, we generalize the notions of a Ramsey scheme to the asymmetric (directed) setting. We begin by recalling some definitions.

\begin{definition}
Let $U$ be a set, and $m \in \mathbb{Z}^{+}$. A \textit{Ramsey scheme} in $m$ colors is a partition of $U \times U$ into $m$ binary relations $\text{Id}, R_{0}, \ldots, R_{m-1}$ such that the following conditions hold:
\begin{enumerate}[label=(\Alph*)]
\item $R_{i}^{-1} = R_{i}$,
\item $R_{i} \circ R_{i} = R_{i}^{c}$, and
\item For all pairs of distinct $0 \leq i, j < m$, $R_{i} \circ R_{j} = \text{Id}^{c}$.
\end{enumerate}

\noindent Here, $\text{Id} = \{ (u, u) : u \in U\}$ is the identity relation over $U$.
\end{definition}

The usual method of constructing the relations $R_{0}, \ldots, R_{m-1}$ is a \textit{guess-and-check} approach due to Comer \cite{ComerMonochrome}, which works as follows. Fix $m \in \mathbb{Z}^{+}$, and let $p \equiv 1 \pmod{2m}$ be a prime. Let $X_{0} := H \leq \mathbb{F}_{p}^{\times}$ be a subgroup of order $(p-1)/m$. Now let $X_{1}, \ldots, X_{m-1}$ be the cosets of $\mathbb{F}_{p}^{\times}/X_{0}$. In particular, as $\mathbb{F}_{p}^{\times}$ is cyclic, we may write $X_{i} = g^{i}X_{0} = \{ g^{am+i} : a \in \mathbb{Z}^{+} \},$ where $g$ is a generator of $\mathbb{F}_{p}^{\times}$. Suppose that $X_{0}, \ldots, X_{m-1}$ satisfy the following conditions:
\begin{enumerate}[label=(\alph*)]
\item $X_{i} = -X_{i}$, for all $0 \leq i \leq m-1$,
\item $X_{i} + X_{i} = \mathbb{F}_{p} \setminus X_{i}$, for all $0 \leq i \leq m-1$, and
\item For all distinct $0 \leq i, j \leq m-1$, $X_{i} + X_{j} = \mathbb{F}_{p} \setminus \{0\}$.
\end{enumerate}

\noindent \\ For each $0 \leq i \leq m-1$, define $R_{i} := \{(x,y) \in \mathbb{F}_{p} \times \mathbb{F}_{p} : x - y \in X_{i} \}.$ Here, the sets $R_{0}, \ldots, R_{m-1}$ are the atoms in our relation algebra. It is easy to check that conditions (a)-(c) on the sets $X_{0}, \ldots, X_{m-1}$ imply that conditions (A)-(C) from the definition of a Ramsey scheme are satisfied for the relations $R_{0}, \ldots, R_{m-1}$. 

We note that condition (b), that $X_{i} + X_{i} = \mathbb{F}_{p} \setminus X_{i}$, indicates that each $X_{i}$ is sum-free. The fact that $p \equiv 1 \pmod{2m}$ implies that $X_{0}$ has even order. It follows that $X_{0}$ is symmetric; i.e., $X_{0} = -X_{0}$. In \cite{Alm2017401AB}, Ramsey schemes were constructed for all $m \leq 2000$ except for $m = 8, 13$, and it was shown that if $p > m^{4} + 5$, then $X_{0}$ contains a solution to $x + y = z.$ In such cases, Comer's construction fails to yield an $m$-color Ramsey  scheme.

A natural generalization is to consider $X_{0}$ to be of odd order, so that $X_{0} \cap -X_{0} = \emptyset.$ (Note that in general, $-X_i = X_{i+m}$.) This yields relations that are anti-symmetric. This motivates the following definition.

\begin{definition} \label{defRamsey}
Let $k$ be odd and $n$ be even, and let $p = nk + 1$ be prime. Let $m = n/2$, and let $g$ be a primitive root modulo $p$. Let $X_{0} = \{ g^{\alpha n} : 0 \leq \alpha < k \}$ be the  subgroup of $\mathbb{F}_{p}^{\times}$ of index $n$. For $1 \leq i < n$, let $X_{i} = \{g^{\alpha n + i} : 0 \leq \alpha < k \}$ be the remaining cosets. Suppose that $X_{0}, \ldots, X_{n-1}$ satisfy the following conditions.
\begin{enumerate}[label=(\alph*)]
\item For all $0 \leq i < n$, $X_{i} + X_{i} = \bigcup_{j \neq i} X_{j}$,

\item For all $0 \leq i < n$, $X_{i} + X_{m+i} = \{0\} \cup \bigcup_{j \not \in \{i, i+m\}} X_{j}$, and

\item For all $0 \leq i < n$ and all $j \not \in \{i, i+m\}$, $X_{i} + X_{j} = \mathbb{F}_{p}^{\times}.$
\end{enumerate}

\noindent Then $\{X_{0}, \ldots, X_{n-1}\}$ forms a \textit{(prime field) directed Ramsey scheme} in $m$ colors.
\end{definition}

\begin{remark}
Conditions (a)-(c) ensure that each coset is sum-free; but otherwise, sumsets are as large as possible. One might think that the exclusions in (b) are unnecessary, but one would be mistaken. They are implied by the exclusion $X_{i} + X_{i} \not \supseteq X_{i}$ (which excludes the transitive relation $(a, b), (b, c), (a,c)$), the \textit{triangle symmetry}
\[
X_{i} + X_{j} \supseteq X_{\ell} \iff X_{j} + X_{\ell + m} \supseteq X_{i},
\]

and the \textit{rotational symmetry}
\[
X_{i} + X_{j} \supseteq X_{\ell} \iff X_{i + a} + X_{j+a} \supseteq X_{\ell + a}.
\]
\end{remark}

We refer to $m$ colors rather than $n$ colors, as the $n$ cosets come in $m$ pairs, where $X_{i} = -X_{i+m}$. One can give a directed $m$-coloring of the edges in $K_{p}$ by labeling the vertices with the elements of $\mathbb{F}_{p}$ with color classes $R_{i} = \{ (x, y) \in \mathbb{F}_{p} \times \mathbb{F}_{p} : x - y \in X_{i} \}$ for $0 \leq i < m$. (So if $X_{0}$ corresponds to the  \textit{directed blue} color label, then $X_{m}$ would be the \textit{backward blue} color label.) Thus, the cosets $X_{0}, \ldots, X_{m-1}$ suffice to determine the directed edge-coloring of $K_{p}$.

In the setting of directed Ramsey schemes, the monochromatic transitive triples (i.e., $(x, y), (y, z), (x, z)$) are forbidden. This suggests a different generalization, where we forbid monochromatic intransitive triples (i.e., $(x, y), (y, z), (z, x)$).

\begin{definition} \label{defAntiRamsey}
Let $k$ be odd, and let $n$ be even, with $p = nk + 1$ prime. Let $m = n/2$. Let $g$ be a primitive root modulo $p$. Let $X_{0} = \{g^{\alpha n} : 0 \leq \alpha < k \}$ be the subgroup of $\mathbb{F}_{p}^{\times}$ of index $n$, and let $X_{1}, \ldots, X_{n-1}$ be the remaining cosets of $\mathbb{F}_{p}^{\times}/X_{0}$, where for some fixed generator $g$ of $\mathbb{F}_{p}^{\times}$, $X_{i} = g^{i}X_{0} = \{ g^{am+i} : a \in \mathbb{Z}^{+} \}.$ Suppose that the following conditions are satisfied:
\begin{enumerate}[label=(\alph*)]
\item For all $0 \leq i < n$, we have that $X_{i} + X_{i} = \bigcup_{j\neq i + m} X_{j}$,

\item For all $0 \leq i < n$, we have that $X_{i} + X_{i+m} = \mathbb{F}_{p}$, and

\item For all $0 \leq i < n$ and all $j \not \in \{i, i+m\}$, we have that $X_{i} + X_{j} = \mathbb{F}_{p}^{\times}.$
\end{enumerate}

\noindent \\ Then we say that the collection $\{X_{0}, \ldots, X_{n-1}\}$ forms a \textit{directed anti-Ramsey scheme} in $m$ colors.
\end{definition}

This leads us to the following question.

\begin{question}[Open]
For which $m > 1$ is there a prime $p$ such that there exists an $m$-color directed Ramsey scheme (respectively, directed anti-Ramsey  scheme) over $\mathbb{F}_{p}$?
\end{question}

\begin{remark}
As a comment for the experts, we refer to the explicit combinatorial structures in Definitions \ref{defRamsey} and \ref{defAntiRamsey} as directed (anti-)Ramsey schemes, and to the associated relation algebras as directed (anti-)Ramsey relation algebras, following Kowalski \cite{KowalskiRamsey} and Maddux \cite{MadduxRamsey1982}. We observe, for instance, that the $1$-color directed anti-Ramsey algebra is the so-called point algebra, and it is not finitely representable.
\end{remark}

We now introduce a key lemma that limits our search space for directed (anti-)Ramsey schemes over a fixed prime $p$.

\begin{lemma} \label{DirectedRamseyBounds}
Let $p = nr + 1$ be prime, with $r$ odd and $n = 2m$. If $p > m^{4} - (2-o(1))m^{3} + (1-o(1))^{2} m^{2} + 4$, then $\mathbb{F}_{p}$ does not yield via Comer's construction an $m$-color directed Ramsey scheme.
\end{lemma}

\begin{proof}
We modify the proof of \cite[Theorem~4]{Alm2017401AB}. We show that if $p > m^{4} - (2-o(1))m^{3} + (1-o(1))^{2}m^{2} + 4$, then $X_{0}$ is not sum-free. Our goal is to count the number $\mathcal{N}$ of solutions to the equation $x + y = z$ inside of $X_{0} \subseteq \mathbb{F}_{p}$. Here, we assume that $|X_{0}| = \delta p$, for $\delta = (p-1)/mp$. By the First Orthogonality Relation for characters, we have that
\[
\frac{1}{p} \sum_{k=0}^{p-1} \chi_{x}(k) = \begin{cases} 1 & : x \equiv 0 \pmod{p}, \\ 0 & : \text{otherwise}. \end{cases}
\]

\noindent Let $\text{Ch}_{X_{0}}$ be the characteristic function of $X_{0}$. It follows that
\begingroup
\allowdisplaybreaks
\begin{align*}
\mathcal{N} &= \frac{1}{p} \sum_{(x,y,z) \in X_{0}^{3}} \sum_{k=0}^{p-1} \chi_{k}(x+y-z) \\
&= \frac{1}{p} \sum_{(x,y,z) \in X_{0}^{3}} \sum_{k=0}^{p-1} \chi_{x+y-z}(k) \\
&= \frac{1}{p} \sum_{(x,y,z) \in X_{0}^{3}} \sum_{k=0}^{p-1} \chi_{x}(k)\chi_{y}(k)\chi_{-z}(k) \\
&= \frac{1}{p} \sum_{(x,y,z) \in X_{0}^{3}} \sum_{k=0}^{p-1} \chi_{x}(k) \chi_{y}(k) \overline{\chi_{z}(k)} \\
&= \frac{1}{p} \sum_{k=0}^{p-1} \left( \sum_{x=0}^{p-1} \text{Ch}_{X_{0}}(k) \chi_{x}(k) \right) \cdot \left( \sum_{y=0}^{p-1} \text{Ch}_{X_{0}}(k) \chi_{y}(k) \right) \cdot \left( \sum_{z=0}^{p-1} \text{Ch}_{X_{0}}(k) \overline{\chi_{z}(k)} \right)\\
&= \frac{1}{p} \sum_{k=0}^{p-1}  \widehat{\text{Ch}_{X_{0}}}(k)^{2} \cdot \widehat{\text{Ch}_{X_{0}}}(-k) \\
&= \frac{1}{p} \cdot \widehat{\text{Ch}_{X_{0}}}(0)^{3} + \frac{1}{p} \cdot  \sum_{k=1}^{p-1}  \widehat{\text{Ch}_{X_{0}}}(k)^{2} \cdot \widehat{\text{Ch}_{X_{0}}}(-k) \\
&= \frac{1}{p} \cdot |X_{0}|^{3} + \frac{1}{p} \cdot  \sum_{k=1}^{p-1}  \widehat{\text{Ch}_{X_{0}}}(k)^{2} \cdot \widehat{\text{Ch}_{X_{0}}}(-k) \\
&= \frac{1}{p} \cdot (\delta p)^{3} + \frac{1}{p} \cdot  \sum_{k=1}^{p-1}  \widehat{\text{Ch}_{X_{0}}}(k)^{2} \cdot \widehat{\text{Ch}_{X_{0}}}(-k) \\
&=\delta^{3} p^{2} + \frac{1}{p} \cdot  \sum_{k=1}^{p-1}  \widehat{\text{Ch}_{X_{0}}}(k)^{2} \cdot \widehat{\text{Ch}_{X_{0}}}(-k) \\
&\leq \delta^{3}p^{2} + \frac{1}{p} \max |\widehat{\text{Ch}_{X_{0}}}(k)| \cdot \left|\sum_{k=1}^{p-1} \widehat{\text{Ch}_{X_{0}}}(k)^{2} \right|,
\end{align*}
\endgroup

\noindent Here, the second equality follows from the fact that $\chi_{a}(k) = \chi_{k}(a)$. We now turn to bounding the term
\begin{align*}
\max| \widehat{\text{Ch}_{X_{0}}}(k)| &\leq \frac{1}{m} \cdot (1 + (m-1)\sqrt{p}) \\
&\leq \frac{1}{m} \cdot (m - (1-o(1))) \sqrt{p}.
\end{align*}

\noindent The first line follows from \cite[Theorem~6.8]{BabaiFourier}, and the second line follows from the fact that $\frac{1}{\sqrt{p}} \leq o(1)$. Thus, we have that
\begin{align*}
\frac{1}{p} \max |\widehat{\text{Ch}_{X_{0}}}(k)| \cdot \left|\sum_{k=1}^{p-1} \widehat{\text{Ch}_{X_{0}}}(k)^{2} \right| 
&\leq \frac{1}{p} \cdot \left( \frac{1}{m} \cdot (m-(1-o(1)))\sqrt{p} \right) \cdot \left|\sum_{k=1}^{p-1} \widehat{\text{Ch}_{X_{0}}}(k)^{2} \right| \\
&\leq \left(\frac{1}{m} \cdot (m-(1-o(1)))\sqrt{p} \right) \cdot  \left|\sum_{k=1}^{p-1} \text{Ch}_{X_{0}}(k)^{2} \right| \\
&\leq \left(\frac{1}{m} \cdot (m-(1-o(1)))\sqrt{p} \right) \delta p \\
\end{align*}

\noindent Here, the second line follows by Parseval's identity. Hence, $\mathcal{N} \geq \delta^{3}p^{2} - \left(\frac{1}{m} \cdot ((m-(1-o(1)))p^{3/2}) \right) \delta$. So there is at least one non-trivial solution to the equation $x + y = z$ in $X_{0}$ if: 
\[
\left(\frac{1}{m} \cdot (m-(1-o(1)))p^{3/2} \right) \delta < \delta^{3}p^{2}. 
\]

\noindent This is equivalent to:
\[
\left(\frac{1}{m} \cdot (m-(1-o(1)))p^{1/2} \right) < \delta^{2}p.
\]

\noindent Now as $\delta = (p-1)/(mp)$, 
\[
\left(\frac{1}{m} \cdot (m-(1-o(1)))p^{1/2} \right) < \delta^{2}p,
\]

\noindent is equivalent to
\[
m(m-(1-o(1)))p^{3/2} < (p-1)^2.
\]

\noindent In particular, we have that
\[
m^2(m-(1-o(1)))^2 p^{3} < (p-1)^{4}.
\]

\noindent So $p > m^{4} - (2-o(1))m^{3} + (1-o(1))^{2} m^{2} + 4$, which holds for $p$ sufficiently large. Therefore, if $p > m^{4} - (2-o(1))m^{3} + (1-o(1))^{2} m^{2} + 4$, then $X_{0}$ contains a solution to $x + y = z$ and hence is not sum-free.
\end{proof}

\begin{lemma}  \label{DirectedAntiRamseyBounds}
Let $p = nr + 1$ be prime, with $r$ odd and $n = 2m$. If $p > m^{4} - (2-o(1))m^{3} + (1-o(1))^{2} m^{2} + 4$, then $\mathbb{F}_{p}$ does not yield via Comer's construction an $m$-color directed anti-Ramsey scheme.
\end{lemma}

\begin{proof}[Proof (Sketch).]
The proof is almost identical to that of Lemma \ref{DirectedRamseyBounds}. The key idea is to replace the equation $x + y = z$ by $x + y = -z$. The argument goes through \textit{mutatis mutandis}.
\end{proof}

\begin{remark}
Lemmas \ref{DirectedRamseyBounds} and \ref{DirectedAntiRamseyBounds} provide that if Comer's construction yields an $m$-color directed (anti-)Ramsey scheme, then $p \leq m^{4} - (2-o(1))m^{3} + (1-o(1))^{2} m^{2} + 4$. This improves upon the previous bounds of $p \leq m^{4} + 5$ from \cite{Alm2017401AB, AlonBourgain}.
\end{remark}

\begin{remark}
In the case when $p > m^{4} - (2-o(1))m^{3} + (1-o(1))^{2} m^{2} + 4$, the (anti-)Ramsey scheme has no forbidden cycles, and so its isomorphism type is determined. 
\end{remark}

In the process of proving Lemma \ref{DirectedRamseyBounds}, we were able to adapt the techniques to improve the upper bound on the multiplicative van der Waerden number introduced in \cite{AlmArithmeticProgressions}. Alm previously showed that $\text{VW}^{\times}(n) \leq (1+o(1))n^{4}$. We improve this bound as follows.

\begin{proposition}
$\text{VW}^{\times}(n) \leq (1-o(1))n^{4} + O(n^{7/2})$.
\end{proposition}

\begin{proof}
For full details, see Appendix \ref{appendix:AlmArithmeticProgressions}.
\end{proof}

\subsection{Data}
\begin{figure}[H]
    \centering
    \includegraphics[width=4in]{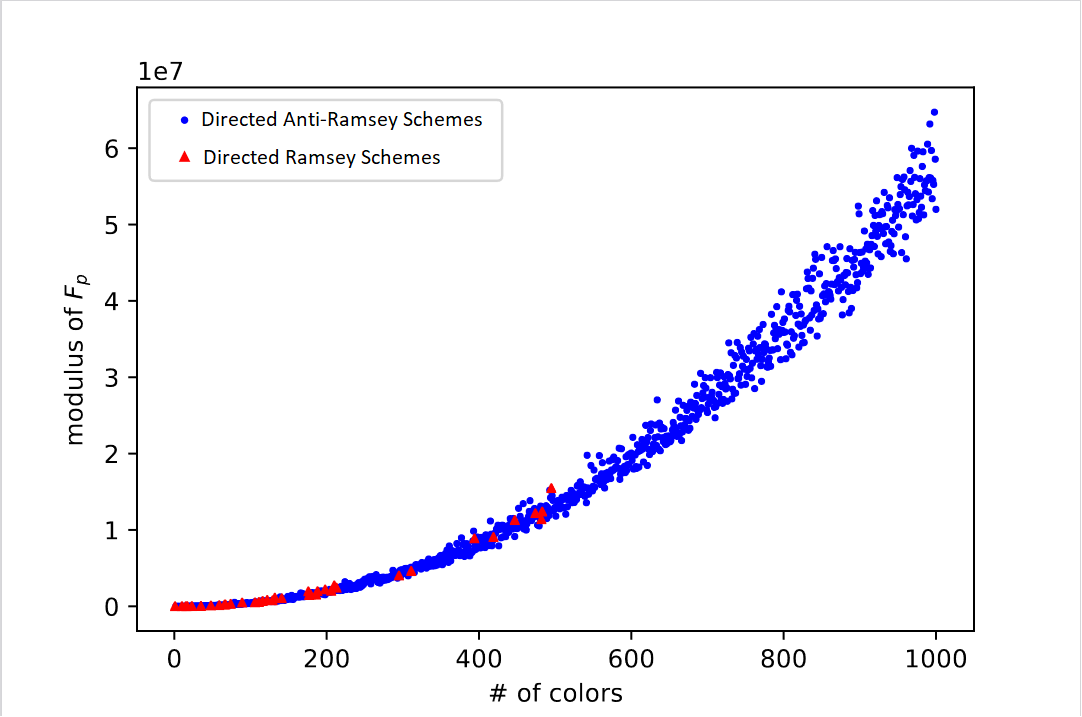}
    \caption{Smallest modulus $p$ for $m$-color directed (anti-)Ramsey schemes}
    \label{fig:moduli}
\end{figure}

We found $m$-color directed Ramsey schemes for 37 values of $m<500$; see Table \ref{tab:my_label}. 
We found $m$-color directed anti-Ramsey schemes for all values of $1<m\leq 1000$ except for $m=8$. All data collected for this paper made use of the fast algorithm for computing cycle structures of Comer's algebras over $\mathbb{F}_{p}$ given in \cite{AlmYlv}. The data for directed anti-Ramsey schemes, as well as the Python code, are available at \href{https://oeis.org/A294615}{https://oeis.org/A294615}. See Figure  \ref{fig:moduli}.

\begin{table}[]
    \centering
   \begin{tabular}{c|c||c|c||c|c}
        $m$ & $p$ & $m$ & $p$ & $m$ & $p$ \\
        \hline
             1 & 3 & 106 & 497141 & 198 & 2192653 \\
        10 & 3221 & 111 & 559219 & 206 & 2020861 \\
        15 & 4231 & 116 & 679993 & 210 & 2728741 \\
        17 & 11527 & 122 & 814717 & 213 & 2420959 \\
        23 & 15319 & 129 & 764971 & 295 & 4017311 \\
        35 & 38011 & 132 & 1118569 & 311 & 4618351 \\
        48 & 91873 & 141 & 1043683 & 394 & 8881549 \\
        59 & 135347 & 176 & 1946209 & 419 & 9071351 \\
        66 & 209221 & 177 & 1470871 & 447 & 11279599 \\
        67 & 228203 & 179 & 1521859 & 474 & 12190333 \\
        74 & 309173 & 186 & 1514413 & 482 & 11383877 \\
        89 & 476863 & 188 & 1968361 & 483 & 12390883 \\
            &       &       &       & 495 & 15468751 \\
    \end{tabular}

    \caption{Smallest prime modulus $p$ for $m$-color directed Ramsey schemes}
    \label{tab:my_label}
\end{table}

The difference in ease in constructing directed Ramsey and anti-Ramsey schemes was at first surprising to the first author. But there is a fundamental asymmetry that is Ramsey-theoretic in nature. For consider the subgraph of the coloring described in Section \ref{sec:RamseySchemes} induced by any one color. With directed anti-Ramsey schemes, this subgraph will consist of multiple transitive triangles pasted together, and can contain large cliques. With directed Ramsey schemes, however, any such color class contains only intransitive triangles. This forces the color class to be sparse in the sense that the \emph{only} cliques are triangles, since any tournament (directed clique) on 4 or more vertices contains a transitive triple $(a,b),(b,c),(a,c)$. This is a plausible explanation for the greater difficulty in finding directed Ramsey schemes.

\section{Finite Relation Algebras} 

We will give some explicit finite representations of relation algebras.  Many of these algebras contain flexible atoms. Thus we eliminate several potential small counterexamples to the Flexible Atom Conjecture.

Since all representations considered here are group representations, throughout we let $\rho$ be a representation that maps into the powerset of $\mathbb{F}_{p}$, where $p$ is given by context.  We use the numbering system and notation given in Maddux's book \cite{Madd}.  Unless otherwise noted, all algebras in this section were not previously known to be finitely representable. In most cases, we use the fact that the algebra in question embeds in some directed Ramsey or anti-Ramsey scheme.

\subsection{The class $a$, $r$, $\Breve{r}$}
\subsubsection{$33_{37}$}

\begin{figure}[H]
    \centering
    \includegraphics[width=4in]{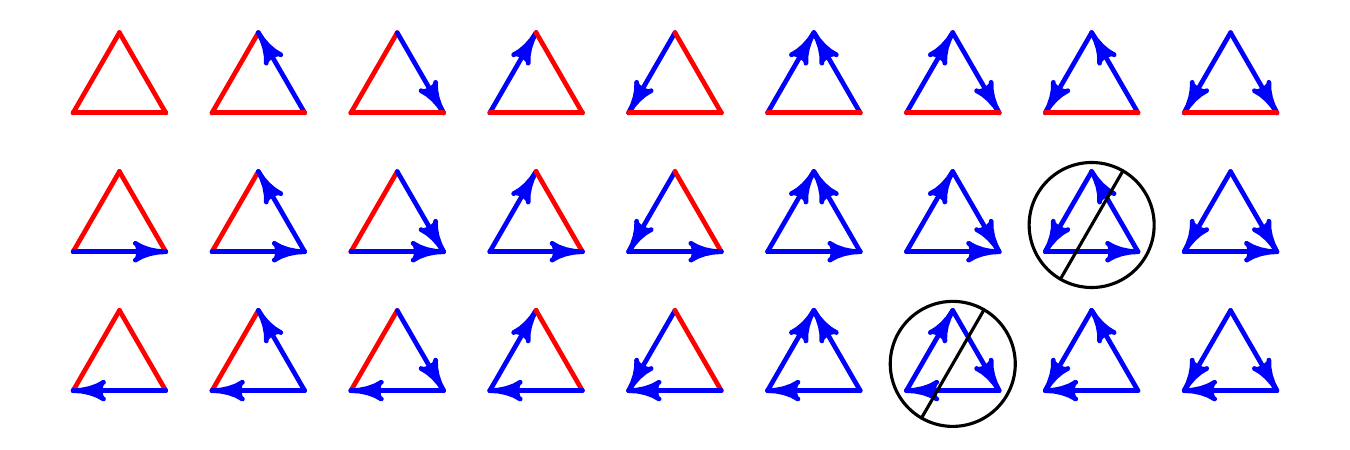}
    \caption{Cycle structure of $33_{37}$, where the entries in a given row outline the needs of the bottom edge and the entries in a column specify relational composition.}
    \label{fig:33_37}
\end{figure}

Relation algebra $33_{37}$ has atoms $1'$, $a$, $r$, $\Breve{r}$. The forbidden diversity cycle is $rr\Breve{r}$. (See Figure \ref{fig:33_37}.)  Let $p=29$ and $m=2$. Then
\begin{align*}
    \rho(a)&=X_1\cup X_3\\
    \rho(r)&=X_0\\
    \rho(\Breve{r})&=X_2
\end{align*} is a representation over $\mathbb{F}_{p}$. Notice that $33_{37}$ is a subalgebra of the 2-color directed anti-Ramsey algebra.

Relation algebra $33_{37}$ was previously shown to be finitely representable by J.~Manske and the first author (unpublished), but the representation given here is smaller. 

\subsubsection{$35_{37}$}

\begin{figure}[H]
    \centering
    \includegraphics[width=4in]{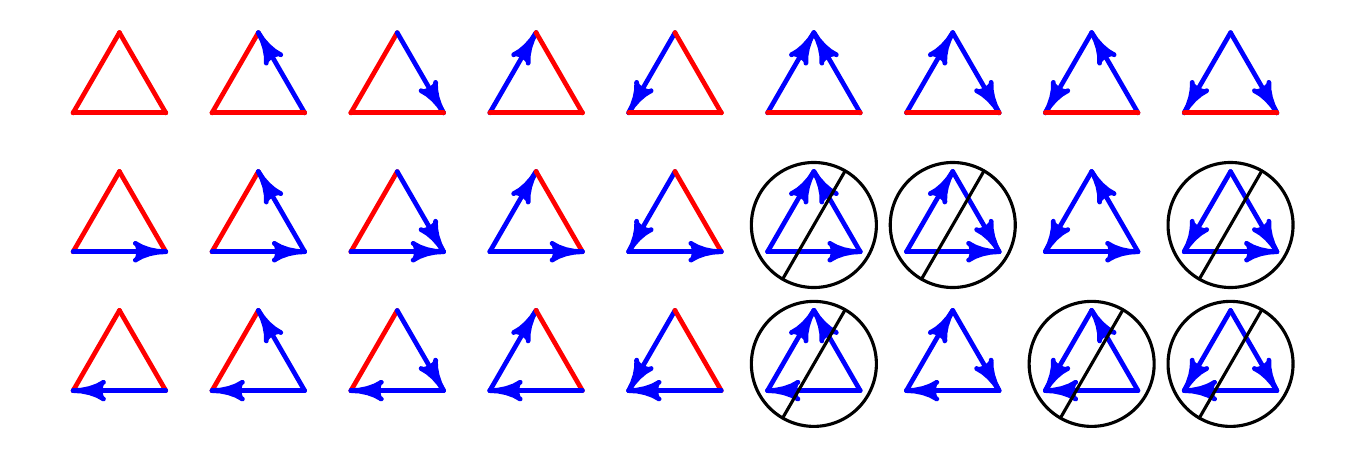}
    \caption{Cycle structure of $35_{37}$, where the entries in a given row outline the needs of the bottom edge and the entries in a column specify relational composition.}
    \label{fig:35_37}
\end{figure}

Relation algebra $35_{37}$ has atoms $1'$, $a$, $r$, $\Breve{r}$. The forbidden diversity cycle is $rrr$. (See Figure \ref{fig:35_37}.)  Let $p=3221$ and $m=10$. Then 
\begin{align*}
    \rho(a)&=\bigcup_{0\neq i \neq 10} X_i\\
    \rho(r)&=X_0\\
    \rho(\Breve{r})&=X_{10}
\end{align*} is a representation over $\mathbb{F}_{p}$. Notice that $35_{37}$ is a subalgebra of the 10-color directed Ramsey algebra.

\subsection{The class  $r$, $\Breve{r}$,  $s$, $\Breve{s}$}

\subsubsection{$77_{83}$}
Relation algebra $77_{83}$ is the 2-color directed anti-Ramsey algebra. The forbidden diversity cycles are $rr\Breve{r}$ and $ss\Breve{s}$. Let $p=29$ and $m=2$. Then 
\begin{align*}
    \rho(r)&=X_0\\
    \rho(\Breve{r})&= X_2\\
    \rho(s)&=X_{1}\\
    \rho(\Breve{s})&=X_{3}
\end{align*} is a representation over $\mathbb{F}_{p}$.

\subsubsection{$78_{83}$}
Relation algebra $78_{83}$ has forbidden diversity cycles  $sss$ and $ss\Breve{s}$. The atoms $r$ and $\Breve{r}$ are flexible.  Let $p=33791$ and $m=31$. Then 
\begin{align*}
    \rho(r)&=\bigcup_{0< i < 31} X_i\\
    \rho(\Breve{r})&=\cup_{31< i < 62} X_i\\
    \rho(s)&=X_{0}\\
    \rho(\Breve{s})&=X_{31}
\end{align*} is a representation over $\mathbb{F}_{p}$. In this case, $\mathbb{F}_{p}$ admits a 31-color (symmetric) Ramsey algebra, but the symmetric colors are splittable into asymmetric pairs. (Because the order of the cosets is even but not divisible by four, each coset may be split into two  cosets of twice the index, and because their order is odd they are asymmetric.)

\subsubsection{$80_{83}$}
Relation algebra $80_{83}$ has forbidden diversity cycle  $ss\Breve{s}$. The atoms $r$ and $\Breve{r}$ are flexible. Let $p=67$ and $m=3$. Then 
\begin{align*}
    \rho(r)&=X_1 \cup X_2\\
    \rho(\Breve{r})&=X_4 \cup X_5\\
    \rho(s)&=X_{0}\\
    \rho(\Breve{s})&=X_{3}
\end{align*} is a representation over $\mathbb{F}_{p}$. $80_{83}$ embeds into the 3-color directed anti-Ramsey algebra.

\subsubsection{$82_{83}$}
Relation algebra $82_{83}$ has forbidden diversity cycle $sss$. The atoms $r$ and $\Breve{r}$ are flexible. Let $p=3221$ and $m=10$. Then 
\begin{align*}
    \rho(r)&=\bigcup_{0 < i < 10} X_i\\
    \rho(\Breve{r})&=\bigcup_{10 < i < 20} X_i\\
    \rho(s)&=X_{0}\\
    \rho(\Breve{s})&=X_{10}
\end{align*} is a representation over $\mathbb{F}_{p}$. $82_{83}$ embeds into the 10-color directed Ramsey algebra.

\subsubsection{$83_{83}$}
Relation algebra $83_{83}$ has no forbidden diversity cycles; all diversity atoms are flexible. $83_{83}$  embeds in the 4-color directed anti-Ramsey algebra; hence, it is representable over $\mathbb{F}_{p}$ for $p=233$.  There is a smaller ``direct'' representation over $\mathbb{F}_{p}$, $p=37$ and $m=2$, where the images of the four diversity atoms are just the four cosets of the multiplicative subgroup of index 4.  This algebra was previously known to be finitely representable, by a slight generalization of the probabilistic argument in \cite{Jipsen}.

\subsection{The class $a$, $b$ $r$, $\Breve{r}$}
\subsubsection{$1310_{1316}$}
Relation algebra $1310_{1316}$ has forbidden diversity cycle $rr\Breve{r}$. The atoms $a$ and $b$ are flexible. Let $p=67$ and $m=3$. Then 
\begin{align*}
    \rho(a) &= X_1 \cup X_4\\
    \rho(b) &= X_2 \cup X_5\\
    \rho(r)&= X_0\\
    \rho(\Breve{r})&= X_3
\end{align*} is a representation over $\mathbb{F}_{p}$. $1310_{1316}$ embeds into the 3-color directed anti-Ramsey algebra.

\subsubsection{$1313_{1316}$}
Relation algebra $1313_{1316}$ has forbidden diversity cycle $rrr$. The atoms $a$ and $b$ are flexible. Let $p=3221$ and $m=10$. Then 
\begin{align*}
    \rho(a) &= \left(\bigcup_{0< i < 6} X_i \right) \cup \left(\bigcup_{10< i < 16} X_i\right)\\
    \rho(b) &= \left(\bigcup_{5< i < 10} X_i \right) \cup \left(\bigcup_{15< i < 20} X_i\right)\\
    \rho(r) &= X_0\\
    \rho(\Breve{r}) &= X_{10}\\
\end{align*} is a representation over $\mathbb{F}_{p}$. $1313_{1316}$ embeds into the 10-color directed Ramsey algebra. 

\subsubsection{$1315_{1316}$} \label{1315}
Relation algebra $1315_{1316}$ has forbidden diversity cycle $bbb$. The atoms $a$, $r$, and $\Breve{r}$ are flexible. Let $p=33791$ and $m=31$. Then 
\begin{align*}
    \rho(a) &= \left(\bigcup_{2< i < 31} X_i \right) \cup \left(\bigcup_{33< i < 62} X_i\right)\\
    \rho(b) &= X_0 \cup X_{31}\\
    \rho(r) &= X_1\cup X_2\\
    \rho(\Breve{r}) &= X_{32}\cup X_{33}\\
\end{align*} is a representation over $\mathbb{F}_{p}$. Note that $\mathbb{F}_p^\times$ has a unique multiplicative subgroup of index 62. Its cosets come in 31 inverse pairs, and these cosets generate a scheme. The underlying relations of this scheme  generate a relation algebra that contains a (symmetric) Ramsey Algebra as a subalgebra. The atoms of this Ramsey algebra behave exactly the way we need the atom $b$ to behave: $b\circ b = \neg b$. So it will suffice to take $\rho(b) = X_0 \cup X_{31}$, since $-X_0 = X_{31}$. To get an antisymmetric flexible atom, we can take $\rho(r) = X_1\cup X_2$, which is ``big enough'' (and  hence $\rho(\Breve{r}) = X_{32}\cup X_{33}$), then $\rho(a)$ is ``everything else."

\subsubsection{$1316_{1316}$} \label{1316}
Relation algebra $1316_{1316}$  has no forbidden diversity cycles; all diversity atoms are flexible. Let $p=73$ and $m=4$. Then 
\begin{align*}
    \rho(a) &= X_2 \cup X_{6}\\
    \rho(b) &= X_3 \cup X_{7}\\
    \rho(r) &= X_0\cup X_1\\
    \rho(\Breve{r}) &= X_{4}\cup X_{5}\\
\end{align*} is a representation over $\mathbb{F}_{p}$.  Note that $\mathbb{F}_p^\times$ has a unique multiplicative subgroup of index 8 whose cosets generate a scheme. The underlying relations of this scheme  generate a relation algebra that contains a subalgebra where all of the diversity atoms are symmetric and flexible. These 8 cosets have odd order, hence are asymmetric.  For atoms $a$ and $b$, we pair a coset with its inverse, as in $\rho(a) = X_2 \cup X_{6}$. For atoms $r$ and $\breve{r}$, we pair a coset with one other than its inverse, as in $\rho(r) = X_0\cup X_1$.

This algebra was previously known to be finitely representable, by a slight generalization of the probabilistic argument in \cite{Jipsen}.

\section{Lower Bounds for $31_{37}, 33_{37},$ and $35_{37}$}

In this section, we establish lower bounds on the minimum sized square representations of $\mathrm{Cy}_{n}, \mathrm{Ch}_{n}$, and $\mathrm{Cych}_{n}$. We note that  $\mathrm{Cy}_{n}, \mathrm{Ch}_{n}$, and $\mathrm{Cych}_{n}$ admit a square representation on $m$ points precisely if we can color an appropriate graph on $m$ vertices that forbids the appropriate substructures. In the cases of $\mathrm{Cy}_{n}, \mathrm{Ch}_{n}$, and $\mathrm{Cych}_{n}$, our graph has $m^{2} - m$ edges.

We color the edges of the graph using red and $n$ shades of blue, $b_{1}, \ldots, b_{n}$. The red edge corresponds to the flexible atom. In cases where we have a red edge between the vertices $u, v$, both $(u, v), (v, u) \in \rho(a)$, where $\rho$ is our representation. For this reason, we draw the red edges as undirected edges.

Now suppose we have an edge $(u, v)$ colored $b_{i}$. Then $(u, v) \in \rho(b_{i})$ and $(v, u) \in \rho(\breve{b_{i}})$. Conversely, if $(u, v) \in \rho(\breve{b_{i}})$, then $(v, u) \in \rho(b_{i})$. For this reason, we only draw one edge labeled $b_{i}$ between $u$ and $v$. This could be the edge $(u, v)$ or the edge $(v, u)$. Both such edges are present, though we only draw one such edge in the graph (this makes it easier to visualize the combinatorial proof for the lower bound). The edge $(u, v)$ indicates that $(u, v) \in \rho(b_{i})$.

Our key technique is to analyze the local combinatorial structure induced by a single red edge. We first introduce a trivial lower bound, which applies to $\mathrm{Cy}_{n}, \mathrm{Ch}_{n}$, and $\mathrm{Cych}_{n}$.

\begin{lemma}  \label{lem:MonoLowerBound}
Let $n \geq 1$. Let $C_{n} \in \{ \mathrm{MonoCy}_{n}, \mathrm{MonoCh}_{n}, \mathrm{MonoCych}_{n}, \mathrm{Cy}_{n}, \mathrm{Ch}_{n}, \mathrm{Cych}_{n}\}$ be a relation algebra. In any representation of $C_{n}$, there are at least $4n^{2} + 4n + 3$ points. 
\end{lemma}

\begin{proof}
Let $x_{0}, x_{1}$ be vertices such that we have the red undirected edge $\{x_{0}, x_{1}\}$. We now count the needs of this red edge. There are $4n^{2}$ blue-blue needs, $4n$ red-blue or blue-red needs, and one red-red need, yielding our lower bound $f(n) \geq 4n^{2} + 4n + 3$. We break this up as follows.
\begin{itemize}
\item Let $\overrightarrow{BB}$ denote the set of vertices $u$ such that $(x_{i}, u)$ is blue for $i = 0, 1$ (so $(x_{i}, u) \in \rho(b)$). We determine $|\overrightarrow{BB}|$. For a given vertex $u$, we have that $(x_{0}, u)$ is colored $b_{j}$ and $(x_{1}, u)$ is colored $b_{k}$, where $j, k \in [n]$. We select $j, k$ independently. So by the rule of product, $|\overrightarrow{BB}| \geq n^{2}$. 

\item Let $\overleftarrow{BB}$ denote the set of vertices $u$ such that $(u, x_{i})$ is blue for $i = 0, 1$ (so $(u, x_{i}) \in \rho(b)$). By similar argument as for $\overrightarrow{BB}$, we have that $|\overleftarrow{BB}| \geq n^{2}$. 

\item Let $\overleftarrow{B}\overrightarrow{B}$ denote the set of vertices $u$ such that $(u, x_{0})$ and $(x_{1}, u)$ are blue (so $(u, x_0), (x_1, u) \in \rho(b)$). By similar argument as for $\overrightarrow{BB}$, we have that $|\overleftarrow{B}\overrightarrow{B}| \geq n^{2}.$

\item Let $\overrightarrow{B}\overleftarrow{B}$ denote the set of vertices $u$ such that $(x_{0}, u)$ and $(u, x_{1})$ are blue (so $(x_{0}, u), (u, x_{1}) \in \rho(b)$). By similar argument as for $\overrightarrow{BB}$, we have that $|\overrightarrow{B}\overleftarrow{B}| \geq n^{2}.$

\item Let $\overrightarrow{(R, B)}$ denote the set of vertices $u$ such that $\{x_{0}, u\}$ is red and $(u, x_{1})$ is blue (so $(u, x_{1}) \in \rho(b)$). Observe that $|\overrightarrow{(R, B)}| \geq n$.

\item Define similarly, $\overrightarrow{(B, R)}$ to be the set of vertices $u$ such that $(x_{0}, u)$ is blue (so $(x_{0}, u) \in \rho(b)$) and $\{x_{1}, u\}$ is red. Observe that $|\overrightarrow{(R, B)}| \geq n$.

\item We define analogously $\overleftarrow{(R, B)}$ to denote the set of vertices $u$ such that $\{x_{0}, u\}$ is red and $(x_{1}, u)$ is blue (so $(x_{1}, u) \in \rho(b)$). Observe that $|\overleftarrow{(R, B)}| \geq n$.

\item We define analogously $\overleftarrow{(B, R)}$ to denote the set of vertices $u$ such that $(u, x_{0})$ is blue (so $(u, x_{0} \in \rho(b)$) and $\{u, x_{1}\}$ is red. Observe that $|\overleftarrow{(B, R)}| \geq n$.
\end{itemize}

\noindent \\ This gives us a lower bound of $f(n) \geq 4n^{2} + 4n + 3$, as desired.
\end{proof}

In the case of $n = 1$, our bound applies to $31_{37}, 33_{37}$, and $35_{37}$. We obtain the following.

\begin{corollary} \label{TrivialLowerBound}
If $k < 11$ and $A \in \{31_{37}, 33_{37}, 35_{37}\}$, then $k \not \in \text{Spec}(A)$.
\end{corollary}

\begin{remark} \label{RemarkGeneralRamsey}
If we forbid both the blue cycle and the blue chain, we can further strengthen this bound. As the classical Ramsey number $R(3, 4) = 9$ and as $f$ is monotone, we have that for all $n \geq 1$, there must exist a red clique on $4$ vertices. Take $x_{0}, x_{1}$ to be two such points in this clique. We iterate on the argument in Lemma~\ref{lem:MonoLowerBound}, modifying it slightly to account for the red-red need. Note that counting the red-blue and blue-blue needs yields $\geq 4n^{2} + 4n$ points. Observe that $\{x_{0},x_{1}\}$ has its $r-r$ need met twice in this red clique. This necessitates two additional points, bringing our total to $4n^{2} + 4n + 2$ points. Now including $x_{0}, x_{1}$, we obtain that  $f(n) \geq 4n^{2} + 4n + 4$. We record this with the following observation.
\end{remark}

\begin{lemma}
In any  representation of $\mathrm{MonoCych}_{n}$ or of $\mathrm{Cych}_{n}$, there are at least $4n^{2} + 4n + 4$ points. 
\end{lemma}

This yields the following corollary.

\begin{corollary}
If $k < 12$, then $k \not \in \text{Spec}(31_{37})$.
\end{corollary}

In the cases of $31_{37}$ and $35_{37}$, we can further improve our lower bounds.

\begin{theorem} \label{thm:3137}
The minimum element of $\text{Spec}(31_{37})$ is at least $14$.
\end{theorem}

\begin{proof}
Let $x_{0}, x_{1}$ be vertices such that $\{x_{0}, x_{1}\}$ is a red edge. We note that $\{x_{0},x_{1}\}$ has $9$ needs, drawn from $\{r, b, \breve{b}\} \times \{r, b, \breve{b}\}$. Including $x_{0}, x_{1}$, this necessitates at least $11$ points. 

Now let $u$ be a point witnessing the $b-b$ need, and let $v$ be a point witnessing the $\breve{b}-\breve{b}$ need. As the blue cycle and the blue chain are forbidden, we have that $\{u, v\}$ is a red edge. Now let $w$ be a point such that for some $i = 0,1$, we have that $(x_{i}, w)$ is colored either $b$ or $\breve{b}$. As the blue cycle and the blue chain are forbidden, $w$ cannot witness a need for $\{u,v\}$ from $\{b, \breve{b}\} \times \{b, \breve{b}\}$. We note that a vertex witnessing the $r-r$ need for $\{x_{0},x_{1}\}$ can witness a need from $\{b, \breve{b}\} \times \{b, \breve{b}\}$. This leaves $3$ unmet needs from $\{b, \breve{b}\} \times \{b, \breve{b}\}$, necessitating at least $3$ additional points. This gives a total of $14$ points.
\end{proof}

\begin{theorem}
The minimum element of $\text{Spec}(35_{37})$ is at least $12$.
\end{theorem}

\begin{proof}
Let $\{x_{0},x_{1}\}$ form a red edge. We note that $\{x_{0},x_{1}\}$ has $9$ needs drawn from $\{r, b, \breve{b}\} \times \{r,b,\breve{b}\}$. Including $x_{0}, x_{1}$, this necessitates $11$ points. 

Now let $u$ be a point witnessing the $b-b$ need for $\{x_{0}, x_{1}\}$, and let $v$ be a point witnessing the $b-r$ need for $\{x_{0},x_{1}\}$. As the blue chain is forbidden, $\{u,v\}$ must form a red edge (otherwise, $(x_{0}, u), (x_{0}, v), (u, v)$ would form a blue chain). So $\{u,v\}$ has $9$ needs. We note that $x_{0}$ witnesses the $\breve{b}-b$ need for $\{u, v\}$, and $x_{1}$ witnesses the $\breve{b}-r$ need for $\{u, v\}$.

We consider the following cases.
\begin{itemize}
\item Suppose that $w$ witnesses the $r-r$ need for $\{x_{0},x_{1}\}$. There is no restriction on the need that $w$ can witness for $\{u,v\}$.

\item Suppose $w$ witnesses the $r-b$ need for $\{x_{0},x_{1}\}$. As the blue chain is forbidden and as $u$ witnesses the $b-b$ need for $\{x_{0},x_{1}\}$, we have that $\{u,w\}$ forms a red edge. So $w$ can witness an $r-r$, $r-b$, or $r-\breve{b}$ need for $\{u,v\}$.

\item Suppose that $w$ witnesses the $r-\breve{b}$ need for $\{x_{0},x_{1}\}$. As the blue chain is forbidden, $w$ can only witness the $\breve{b}-r, \breve{b}-b$, and $\breve{b}-\breve{b}$ needs for $\{u,v\}$.

\item Suppose that $w$ witnesses the $\breve{b}-r$ need for $\{x_{0},x_{1}\}$. As the blue chain is forbidden, $w$ can only witness the $b-b, b-r, r-b, r-r$ needs for $\{u,v\}$.

\item Suppose that $w$ witnesses the $\breve{b}-b$ need for $\{x_{0},x_{1}\}$. As the blue chain is forbidden, $\{u,w\}$ necessarily forms a red edge. Furthermore, we cannot have a blue $(w, v)$ edge. So $w$ can only witness the $r-\breve{b}$ and $r-r$ needs for $\{u,v\}$.

\item Suppose that $w$ witnesses the $b-\breve{b}$ need for $\{x_{0},x_{1}\}$. As the blue chain is forbidden, $w$ can only witness the $r-r$ need for $\{u,v\}$.

\item Suppose that $w$ witnesses the $\breve{b}-\breve{b}$ need for $\{x_{0},x_{1}\}$. As the blue chain is forbidden, neither $(w, u)$ nor $(w, v)$ can be blue. So $w$ can witness only the $r-r, b-r, r-b, b-b$ needs for $\{u,v\}$.
\end{itemize}

\noindent Observe that the $b-\breve{b}$ need cannot be met with the points we have considered so far. Thus, we must have at least $1$ additional point, bringing our total to $12$ points. 
\end{proof}

\subsection{Lower Bound for $\mathrm{Cych}_{n}$}

In this section, we consider the family of relation algebras $(\text{Cych}_{n})_{n \geq 2}$. In this case, we are able to adapt the techniques from \cite{alm2021improved} to further improve the lower bound. We denote $f_{\mathrm{Cych}}(n) := \min(\Spec(\mathrm{Cych}_{n}))$.

\begin{observation} \label{LemmaForbiddenBlue}
Let $n \geq 2$. Consider the relation algebra $\mathrm{Cych}_{n}$. Let $\{x_{0}, x_{1}\}$ form a red edge. Fix $i = 0, 1$. Suppose that $w, z$ are points such that $(x_{i}, w)$ and $(x_{i}, z)$ are both blue. Then $(w, z)$ and $(z, w)$ are both colored red. 
\end{observation}

\begin{proof}
As the blue cycle and the blue chain are forbidden, neither $(w, z)$ nor $(z, w)$ can be blue. The result follows.
\end{proof}

\begin{theorem} \label{ImprovedCyChLowerBound1}
Let $n \geq 2$. We have that $f_{\text{Cych}}(n) \geq 8n^{2} + 8n + 3$. 
\end{theorem}

\begin{proof}
Label the shades of blue $b_{1}, \ldots, b_{n}$. Now let $u \in \overrightarrow{BB}$ such that $(x_{0}, u), (x_{1}, u)$ are both colored $b_{1}$. Let $v \in \overrightarrow{BB}$ such that $(x_{0}, v), (x_{1}, v)$ are both colored $b_{2}$. As a blue chain is forbidden, we have necessarily that $\{u, v\}$ is a red edge. Thus, $\{u,v\}$ has $(2n+1)^{2}$ needs. We note that $x_{0}$ witnesses the $\overleftarrow{b_{1}} - \overleftarrow{b_{1}}$ need for $\{u,v\}$, while $x_{1}$ witnesses the $\overleftarrow{b_{2}} - \overleftarrow{b_{2}}$ need for $\{u,v\}$. 

As $u, v \in \overrightarrow{BB}$, we have by Observation~\ref{LemmaForbiddenBlue} that for any vertex $w \in \overrightarrow{BB}$ (where $w \not \in \{u,v\}$), that $(u, w), (w, u),$ $(v, w), (w, v)$ must all be red. Similarly, as the blue chain and blue cycle are forbidden, we have that for any vertex $w \in \overleftarrow{BB} \cup \overleftarrow{B}\overrightarrow{B} \cup \overrightarrow{B}\overleftarrow{B}$, $\{u, w\}, \{v, w\}$ must both be red. In particular, $\overrightarrow{BB} \cup \overleftarrow{BB} \cup \overleftarrow{B}\overrightarrow{B} \cup \overrightarrow{B}\overleftarrow{B}$ induces a red clique on $4n^{2}$ points.

Similarly, by Observation~\ref{LemmaForbiddenBlue}, we again have that any vertex in $\overleftarrow{(R, B)}$, $\overleftarrow{(B, R)}$, $\overrightarrow{(R, B)}$, and $\overrightarrow{(B, R)}$ must witness red-red needs for $uv$. 

Thus, there must be at least $(2n+1)^{2} - 2$ additional points. This gives a total of $4n^{2} + 4n + 4 + (2n+1)^{2} - 2 = 8n^{2} + 8n + 3$ points.
\end{proof}

\begin{remark}
The proof of Theorem \ref{ImprovedCyChLowerBound1}  relies crucially on the assumption that $n \geq 2$ to obtain two vertices $u, v \in \overrightarrow{BB}$. In the case of $n = 1$, which is the case of $31_{37}$, we are unable to deduce that $\overrightarrow{BB}$ has at least two elements. 
\end{remark}

\subsection{Lower Bound for $\mathrm{Ch}_{n}$}

In this section, we consider the family of relation algebras $(\text{Ch}_{n})_{n \geq 2}$. We are again able to adapt the techniques from \cite{alm2021improved} to further improve the lower bound. Denote $f_{\mathrm{Cych}}(n) := \min(\Spec(\mathrm{Cych}_{n}))$.

\begin{theorem} 
Let $n \geq 2$. We have that $f_{\mathrm{Ch}}(n) \geq 6n^{2} + 8n + 5$. 
\end{theorem}

\begin{proof}
Let $\{x_{0}, x_{1}\}$ form a red edge. Let $u \in \overrightarrow{BB}$ such that $(x_{i}, u)$ is colored $b_{1}$ for $i = 0, 1$. Let $v \in \overrightarrow{BB}$ such that $(x_{i}, v)$ is colored $b_{2}$ for $i = 0, 1$. As the blue chain is forbidden, we have that $\{u,v\}$ forms a red edge. So $\{u,v\}$ has $(2n+1)^{2}$ needs. We note that $x_{0}$ witnesses the $\overleftarrow{b_{1}}-\overleftarrow{b_{1}}$ need, and $x_{1}$ witnesses the $\overleftarrow{b_{2}}-\overleftarrow{b_{2}}$ need. 

We consider the following cases.
\begin{itemize}
\item \textbf{Case 1:} By similar reasoning as to why $\{u,v\}$ forms a red edge, we have that any vertex in $\overrightarrow{BB}$ witnesses a $r-r$ need for $\{u,v\}$.

\item \textbf{Case 2:} Suppose $w \in \overleftarrow{BB}$. As the blue chain is forbidden, neither $(w, u)$ nor $(w, v)$ can be colored blue. So $w$ can witness a $r-r$, $r-\overleftarrow{b}$, $\overleftarrow{b}-r$, or $\overleftarrow{b}-\overleftarrow{b}$ needs for $\{u,v\}$. This yields a total of $(n+1)^{2} = n^{2} + 2n + 1$ possible needs. However, $|\overleftarrow{BB}| = n^{2}$. So only $n^{2}$ of these needs can be satisfied.

\item \textbf{Case 3:} Suppose $w \in \overleftarrow{B}\overrightarrow{B}$. As the blue chain is forbidden, $(w, u)$ and $(w, v)$ cannot be blue. Similarly, $(v, w)$ and $(v, u)$ cannot be blue. Thus, $w$ must witness a $r-r$ need for $\{u,v\}$.

\item \textbf{Case 4:} Suppose $w \in \overrightarrow{B}\overleftarrow{B}$. By similar argument as in Case 3, $w$ must witness a $r-r$ need for $\{u,v\}$.

\item \textbf{Case 5:} Suppose $w \in \overrightarrow{RB} \cup \overrightarrow{BR}$. As the blue chain is forbidden, $(w, u), (u, w), (w, v), (v, w)$ cannot be blue. So $w$ must witness a $r-r$ need for $\{u,v\}$.

\item \textbf{Case 6:} Suppose $w \in \overleftarrow{RB} \cup \overleftarrow{BR}$. As the blue chain is forbidden, $(w, u)$ and $(w, v)$ cannot be blue. So $w$ can witness the following needs for $\{u,v\}$: $\overleftarrow{b}-\overleftarrow{b}$, $\overleftarrow{b}-r$, $r-\overleftarrow{b}$, $r-r$. This yields a total of $(n+1)^{2} = n^{2} + 2n + 1$ possible needs. However, $|\overleftarrow{RB} \cup \overleftarrow{BR}| = 2n$. So only $2n$ of the $(n+1)^{2}$ possible needs can be met.
\end{itemize}

\noindent \\ We note in particular that none of the points we have considered thus far can satisfy needs for $\{u,v\}$ of the form $\overrightarrow{b}-r$, $r-\overrightarrow{b}$, or $\overrightarrow{b}-\overrightarrow{b}$. So we require at least $n^{2} + 2n$ additional points, for a total of 
\begin{align*}
(2n+1)^{2} + 2 + 2n + 1 + n^{2} + 1 + n^{2} + 2n &= (2n+1)^{2} + 2n^{2} + 4n + 4 \\
&= 6n^{2} + 8n + 5
\end{align*}

\noindent points.
\end{proof}

\section{Discussion and Conclusion}

We gave finite representations of several relation algebras, including $33_{37}, 35_{37}, 77_{83}, 78_{83}, 80_{83}$, $82_{83}, 83_{83}$, $1310_{1316}$, $1313_{1316}, 1315_{1316}$, and $1316_{1316}$. Prior to our paper, only $83_{83}$ and $1316_{1316}$ were known to be finitely representable, due to a slight generalization of \cite{Jipsen}. 

In order to construct our finite representations, we generalized the notion of a Ramsey scheme from \cite{ComerCombinatorial, KowalskiRamsey} to the directed (anti-symmetric) setting. In the process, we established a necessary condition for Comer's construction to yield an $m$-color directed (anti-)Ramsey scheme over $\mathbb{F}_{p}$-- namely, $p < m^{4} - (2-o(1))m^{3} + (1-o(1))^{2}m^{2} + 4$. This improved the previous bound of $p < m^{4} + 5$ \cite{Alm2017401AB, AlonBourgain}. In the process, we established a new bound on certain van der Waerden-like numbers $\text{VW}^{\times}(n) \leq (1-o(1))n^{4} + O(n^{7/2})$, which improved upon the previous bound of $(1+o(1))n^{4}$ due to Alm \cite{AlmArithmeticProgressions}.

We next extended the techniques of \cite{alm2021improved} to improve the lower bounds of two different generalizations of $31_{37}, 33_{37}$, and $35_{37}$. In the process, we obtained that any square representation of $31_{37}$ requires at least $16$ points, any square representation of $33_{37}$ requires at least $11$ points, and any square representation of $35_{37}$ requires at least $14$ points. Furthermore, we provide lower bounds on relation algebras generalizing $31_{37}, 33_{37}, 35_{37}$ in both the directions of Ramsey schemes and \cite{AMM}.

\noindent We conclude with several open problems.

In \Lem{DirectedRamseyBounds}, we were able to improve the upper bound from $p < m^{4} + 5$ to $p < m^{4} - (2-o(1))m^{3} + (1-o(1))^{2}m^{2} + 4$ by squeezing out better bounds on the maximum modulus for the second Fourier coefficient. It seems unlikely that $p \in O(m^{4})$ is the correct upper bound, and we conjecture that $p \in O(m^{3}).$ As a starting point, we ask the following.

\begin{problem}
Improve the upper bound on \Lem{DirectedRamseyBounds} to $p < m^{4-o(1)}$.
\end{problem}

Such an improvement is also likely to yield a better upper bound on $\text{VW}^{\times}(n)$; see \Prop{prop:MultiplicativeVWN}.

It would also be of interest to establish whether directed (anti-)Ramsey schemes exist for all $m$ sufficiently large. Current work in this direction involves computer-aided searches for constructions. While such computation is helpful in constructing examples, it is not sufficient to prove the existence of infinitely many constructions. To this end, we ask the following.

\begin{problem}
Give an infinite family of pairs $(m, p)$, where $p$ is prime, $m$ is a divisor of $(p-1)$, and there exists an $m$-color directed (anti-)Ramsey scheme over $\mathbb{F}_{p}.$ 
\end{problem}

We showed that the minimum-sized representation of $33_{37}$ is in $\{11, \ldots, 29\}$. As the upper and lower bounds are quite close, we believe that the following problem is within reach. 

\begin{problem}
Determine the minimum-sized representation of $33_{37}$.
\end{problem}

Improving either the upper or lower bounds on the minimum element in $\text{Spec}(33_{37})$ would also be of interest. More generally, it would be of interest to determine upper bounds on the minimum-sized representation for $\mathrm{MonoCy}_{n}, \mathrm{MonoCh}_{n}, \mathrm{MonoCych}_{n}, \mathrm{Cy}_{n}, \mathrm{Ch}_{n}$, and $\mathrm{Cych}_{n}$. This brings us to our next problem.

\begin{problem}
Determine whether any of the following families of relation algebras admit polynomial-sized representations: $\mathrm{MonoCy}_{n}, \mathrm{MonoCh}_{n}, \mathrm{MonoCych}_{n}, \mathrm{Cy}_{n}, \mathrm{Ch}_{n}$, and $\mathrm{Cych}_{n}$.
\end{problem}

\appendix
\section{Arithmetic Progressions of Length $3$ in $\mathbb{F}_{p}$} \label{appendix:AlmArithmeticProgressions}

In this section, we consider arithmetic progressions of length $3$ over $\mathbb{F}_{p}$. Alm \cite{AlmArithmeticProgressions} previously introduced the notion of a multiplicative van der Waerden number, which we recall below.

\begin{definition}[{\cite[Definition~1]{AlmArithmeticProgressions}}]
Let $\text{VW}^{\times}(n)$ denote the least prime $q \equiv 1 \pmod{n}$ such that for all primes $p \geq q$ satisfying $p \equiv 1 \pmod{n}$, the multiplicative subgroup $H$ of $\mathbb{F}_{p}^{\times}$ of index $n$ contains a mod-$p$ arithmetic progression (that is, a triple $x, y, z \in H$ such that $x + y = 2z \in H$).
\end{definition}

Chang \cite{Chang2016ArithmeticPI} showed that if $H < \mathbb{F}_{p}^{\times}$ and $|H| > cp^{3/4}$ for some constant $c$, then $H$ contains non-trivial arithmetic $3$-progressions (that is, arithmetic progressions $(a, a+b, a+2b)$, where $b \not \equiv 0 \pmod{p}$). This implies that $\text{VW}^{\times}(n) \in O(n^{4})$. We note that Chang's proof does not make the constant explicit. Alm \cite{AlmArithmeticProgressions} gave a more explicit constant\footnote{Alm \cite{AlmArithmeticProgressions} previously stated that Chang's \cite{Chang2016ArithmeticPI} result implied his upper bound of $(1+\epsilon)n^{4}$. It is not clear that Chang's result implies Alm's. Rather, Chang's result only appears to provide an asymptotic upper bound of $O(n^{4})$.} of $\text{VW}^{\times}(n) \leq (1+o(1))n^{4}$.

We improve the bound on $\text{VW}^{\times}(n)$.

\begin{proposition} \label{prop:MultiplicativeVWN}
$\text{VW}^{\times}(n) \leq (1-o(1))n^{4} + O(n^{7/2})$.
\end{proposition}

\begin{proof}
The proof is similar to \cite[Theorem~2]{AlmArithmeticProgressions}. Let $\delta = (p-1)/(pn)$, and let $A \subseteq \mathbb{F}_{p}$ be a set of size $\delta p$. Let $\mathcal{N}$ be the number of (possibly trivial) solutions to $x + y = 2z$ inside of $A$. By the proof of \cite[Theorem~2]{AlmArithmeticProgressions}, we have that:
\begin{align*}
\mathcal{N} &= \frac{|A|^{3}}{3} + \frac{1}{p} \sum_{k=1}^{p-1} \widehat{\text{Ch}_{A}}(k)^{2} \cdot \widehat{\text{Ch}_{A}}(-2k) \\
&= \delta^{3}p^{2} + \frac{1}{p} \sum_{k=1}^{p-1} \widehat{\text{Ch}_{A}}(k)^{2} \cdot \widehat{\text{Ch}_{A}}(-2k)
\end{align*}

\noindent We now turn to bounding:
\begin{align*}
\left|\frac{1}{p} \sum_{k=1}^{p-1} \widehat{\text{Ch}_{A}}(k)^{2} \cdot \widehat{\text{Ch}_{A}}(-2k) \right| 
\end{align*}

\noindent Here, we now assume that $A$ is a multiplicative subgroup of $\mathbb{F}_{p}^{\times}$. So we have that:
\begin{align*}
\max|\widehat{\text{Ch}_{A}}(k)| &\leq \frac{1}{n}(1 + (n-1)\sqrt{p}) \\
&\leq \left( \frac{n-(1-o(1))}{n} \right)\sqrt{p} \\
&= (1-o(1))\sqrt{p},
\end{align*}

\noindent where the first inequality follows from the proof of \cite[Theorem~6.8]{BabaiFourier}. It now follows that:

\begin{align*}
\left|\frac{1}{p} \sum_{k=1}^{p-1} \widehat{\text{Ch}_{A}}(k)^{2} \cdot \widehat{\text{Ch}_{A}}(-2k) \right| &\leq \frac{1}{p} \cdot \max|\widehat{\text{Ch}_{A}}(k)| \cdot \left|\sum_{k=1}^{p-1} \widehat{\text{Ch}_{A}}(k)^{2} \right| \\
&\leq \frac{1}{p} (1-o(1))\sqrt{p} \cdot \left|\sum_{k=1}^{p-1} \widehat{\text{Ch}_{A}}(k)^{2} \right| \\
&\leq (1-o(1))p^{1/2} \cdot \left|\sum_{k=1}^{p-1} \text{Ch}_{A}(k)^{2} \right|\\
&\leq \delta(1-o(1))p^{3/2}.
\end{align*}

\noindent The second-to-last line follows by Parseval's identity. Thus, $\mathcal{N} \geq \delta^{3}p^{2} - \delta(1-o(1))p^{3/2}$. Now there are $\delta p$ trivial solutions. So subtracting those out, we obtain:
\[
\mathcal{N} - \delta p \geq \delta^{3}p^{2} - \delta p - \delta(1-o(1))p^{3/2}.
\]

\noindent So there is at least one non-trivial solution if:
\[
\delta^{3}p^{2} > \delta p + \delta(1-o(1))p^{3/2}.
\]

\noindent This is equivalent to:
\[
\delta^{2}p > 1 + (1-o(1))p^{1/2}.
\]

\noindent As $\delta = (p-1)/(np)$, we have that:
\begin{align*}
(p-1)^{2} > n^{2}p(1 + (1-o(1))p^{1/2}).
\end{align*}

\noindent Equivocally, we have that:
\begin{align*}
(p-1)^{2} > n^{2}p + (1-o(1))n^{2}p^{3/2}. 
\end{align*}

\noindent Squaring both sides, we obtain:
\[
(p-1)^{4} > n^{4}p^{2} + 2(1-o(1))n^{4}p^{5/2} + (1-o(1))n^{4}p^{3}.
\]

\noindent Thus:
\[
p > (1-o(1))n^{4} + n^{4}/p + (2-o(1))n^{4}p^{-1/2} + O(1).
\]

\noindent As $n$ is a proper divisor of $p-1$, we have that:
\[
n^{4}/p + (2-o(1))n^{4}p^{-1/2} + O(1) \in O(n^{7/2}).
\]

\noindent Thus, if:
\[
p > (1-o(1))n^{4} + O(n^{7/2}),
\]

\noindent then $A$ has a non-trivial mod-$p$ arithmetic progression of length $3$. The result now follows.
\end{proof}

\begin{remark}
As discussed in \cite{AlmArithmeticProgressions}, it appears that $\text{VW}^{\times}(n) \in o(n^{3})$. Improving the upper bound remains of interest. Even showing that $\text{VW}^{\times}(n) \leq n^{4-o(1)}$ would be a significant step in this direction.
\end{remark}

\bibliographystyle{alphaurl}
\bibliography{references}

\end{document}